\documentclass{amsart}
\usepackage{amssymb}
\usepackage{color,graphics,epic}
\usepackage[matrix,arrow]{xy}

\numberwithin{equation}{section}

\theoremstyle{plain}
\newtheorem{thm}{Theorem}[section]

\newtheorem{theorem}[thm]{Theorem}

\newtheorem{lemma}[thm]{Lemma}

\newtheorem{conjecture}[thm]{Conjecture}

\newtheorem{remark}[thm]{Remark}

\newtheorem{example}[thm]{Example}

\newcommand{\lb}{[[}
\newcommand{\rb}{]]}
\newcommand{\R}{\mathbb{R}}
\newcommand{\Z}{\mathbb{Z}}

\newcommand{\C}{\mathbb{C}}

\newcommand{\half}{{\textstyle\frac{1}{2}}}

\newcommand{\htp}{\simeq}
\newcommand{\iso}{\cong}

\newcommand{\A}{\mathcal{A}}
\newcommand{\CF}{\mathit{CF}}

\newcommand{\F}{\mathcal{F}}

\newcommand{\Aff}{\mathcal{A}\mathit{ff}}
\newcommand{\PP}{\mathcal{P}}
\newcommand{\Cont}{\mathcal{C}}
\newcommand{\Olo}{\mathcal{O}}
\newcommand{\Kan}{\mathcal{K}}

\newcommand{\val}{\mathrm{val}}
\newcommand{\Der}{\mathcal{D}}

\newcommand{\End}{\mathit{End}}
\newcommand{\Proj}{\mathrm{Proj}}

\newcommand{\XX}{\mathcal{X}}
\newcommand{\RR}{\mathcal{R}}

\begin{document}
\title[Pairs-of-pants decompositions]{Some speculations on pairs-of-pants\\ decompositions and Fukaya categories}
\author{Paul Seidel}
\date{April 2010}
\maketitle

\section{Introduction\label{sec:intro}}

This is partly a survey and partly a speculative article, concerning a particular question about Fukaya categories of symplectic manifolds. Namely, can we decompose a symplectic manifold into standard pieces, and then reconstruct its Fukaya category by gluing together categories depending only on the geometry of each piece, in a (loosely understood) sheaf-theoretic way? For this to work, some degree of control over pseudo-holomorphic curves is required, an issue which depends on the geometry of the decomposition under consideration. Sheaf-theoretic ideas have been successfully applied to the symplectic geometry of cotangent bundles, starting with the work of Fukaya-Oh \cite{fukaya-oh98}, and followed by Kasturirangan-Oh \cite{oh97c, kasturirangan-oh00} and Nadler-Zaslow \cite{nadler-zaslow06, nadler06} (for a survey of the last-mentioned work and related ideas of Fukaya-Smith, see \cite{fukaya-seidel-smith07b}). Recently, Kontsevich \cite{kontsevich-talk} has proposed a generalization to Stein manifolds whose Lagrangian skeleta have certain singularities. However, that is not quite the direction we wish to take here.

Instead, we focus on situations inspired by homological mirror symmetry \cite{kontsevich94}. In its most naive form, the construction of the mirror via SYZ duality \cite{strominger-yau-zaslow96, leung05} is local. This leads to a sheaf-theoretic description of Fukaya categories of Lagrangian torus fibrations, proposed and partially proved by Kontsevich-Soibelman \cite{kontsevich-soibelman00}; a closely related approach was pursued by Fukaya \cite{fukaya02c}. However, whenever the SYZ fibration has singularities, instanton corrections emanate from these and spread over the base \cite{kontsevich-soibelman06b, gross-siebert06b}, modifying the naive construction by non-local contributions. Our specific aim here is to look at some other situations where instanton corrections may be absent (as the word {\em may} indicates, this is strictly conjectural). In Section \ref{sec:zariski} we consider two instances, namely hypersurfaces in $(\C^*)^n$ and in an abelian variety. In both cases, the pieces of our decomposition are higher-dimensional pairs-of-pants \cite{mikhalkin04}, and the mirror description is in terms of categories of Landau-Ginzburg branes \cite{orlov04} (see also \cite{eisenbud80,buchweitz86}). In Section \ref{sec:analytic} we make an attempt to merge this with \cite{kontsevich-soibelman00}, which means passing from algebraic to nonarchimedean analytic geometry. Actually, we will do so only in the lowest-dimensional case, so the outcome is a new conjectural model for the Fukaya categories of closed surfaces. This is obviously strongly influenced by the proposed homological mirror symmetry for surfaces of genus $\geq 2$, due to Kapustin-Katzarkov-Orlov-Yotov \cite{kapustin-katzarkov-orlov-yotov09} (some cases have now been proved \cite{seidel08b, efimov09}).

{\em Acknowledgments.} This paper owes a somewhat larger than usual debt to others: Mohammed Abouzaid, Denis Auroux, Ludmil Katzarkov, and Nick Sheridan. In particular, Sheridan explained his unpublished work on signed affine structures \cite{sheridan10} to me, which significantly influenced the material in Section \ref{sec:analytic}. NSF provided partial financial support through grant DMS-0652620.

\section{The Zariski topology\label{sec:zariski}}

\subsection*{Hypersurfaces in the complex torus}

We start by stating a version of the homological mirror symmetry conjecture for hypersurfaces in the complex torus $(\C^*)^n$. This amounts to using Mikhalkin's pair-of-pants decompositions \cite{mikhalkin04} and the closely related theory of toric degenerations, for which see for instance \cite{nishinou-siebert06}. Besides that, it is strongly influenced by other versions of mirror symmetry considered by Gross-Siebert \cite{gross-siebert03, gross-siebert06} and Abouzaid \cite{abouzaid05,abouzaid09b}.

Suppose that we have a function $\phi_\Z: \Z^n \rightarrow \Z \cup \{\infty\}$ which is finite at a finite nonempty set of points. To this and a choice of complex parameter $0 < |\epsilon| < 1$ one associates the Laurent polynomial
\begin{equation}
F(z_1,\dots,z_n) = \sum_{\phi_\Z(v) < \infty} \epsilon^{\phi_\Z(v)} z_1^{v_1}\cdots z_{n}^{v_{n}},
\end{equation}
which defines a hypersurface $M = F^{-1}(0) \subset (\C^*)^n$. We will now make some further assumptions. First of all, the set of points where $\phi_\Z < \infty$ should be of the form $P_\Z = P \cap \Z^n$, where $P$ is a bounded integer polytope with nonempty interior. Moreover, there should be a decomposition of $P$ into simplices $P_k$ which are integral and of volume $1/n!$, with the following properties. Extend $\phi_\Z$ to a function $\phi: P \rightarrow \R$ which is affine on each $P_k$. We ask that $\phi$ should be convex, and not differentiable at any point which lies on more than one $P_k$. Then, if we take $\epsilon$ small, $M$ is smooth, and comes with a decomposition into $(2n-2)$-dimensional pairs-of-pants, one for each $P_k$ \cite[Theorem 4]{mikhalkin04}.

The Legendre transform of $\phi$ is the function $\psi: \R^n \rightarrow \R$ defined by
\begin{equation} \label{eq:legendre}
\psi(v) = \max_{w \in P} v \cdot w - \phi(w) = \max_{w \in P_\Z} v \cdot w - \phi_\Z(w).
\end{equation}
This is again continuous and convex. There is a decomposition of $\R^n$ into finitely many integer polytopes $Q_k$, of which some will be noncompact, such that $\psi$ is affine on each $Q_k$, and not differentiable at any point lying on more than one $Q_k$. One can use $\psi$ to define a toric degeneration, whose generic fibres are isomorphic to $(\C^*)^{n}$, as follows. For any integer $k \geq 1$, let $\RR^k$ be the complex vector space whose basis elements correspond to points $v \in \frac{1}{k}\Z^{n+1}$ such that $v_{n+1} \geq \psi(v_1,\dots,v_{n})$. These spaces form a commutative graded ring $\RR$, where the product of generators $v \in \frac{1}{k} \Z^{n+1}$, $w \in \frac{1}{l} \Z^{n+1}$ is
\begin{equation} \label{eq:product}
\frac{k v + l w}{k + l} \in \textstyle\frac{1}{k+l}\Z^{n+1}.
\end{equation}
$\RR$ is also naturally a free module over $\C[t]$, where $t$ acts by mapping $v \in \frac{1}{k} \Z^{n+1}$ to $(v_1,\dots,v_{n},v_{n+1} + \frac{1}{k})$. Hence,
\begin{equation} \label{eq:proj}
\XX = \Proj(\RR)
\end{equation}
is a toric variety which comes with a flat morphism $W: \XX \rightarrow \mathbb{A}^1 = \mathrm{Spec} (\C[t])$. Among the symmetries of that variety, we want to single out the $\C^*$-action for which the monomial corresponding to $v \in \frac{1}{k}\Z^{n+1}$ has weight $k v_{n+1}$. By definition, $W$ is homogeneous with respect to this action. While these properties would be true in greater generality, our assumptions ensure that $\XX$ is smooth, and covered by toric Zariski charts $\C^{n+1} \hookrightarrow \XX$, one for each simplex in the original decomposition, such that $W(z) = -z_1 \cdots z_{n+1}$ in each chart (however, the weights of the $\C^*$-action can vary from one chart to the next).

\begin{example}
For $n = 2$, $M$ is a curve whose genus is the number of integer points in the interior of $P$. For instance, take the following $\phi_\Z$ and its associated triangulation:
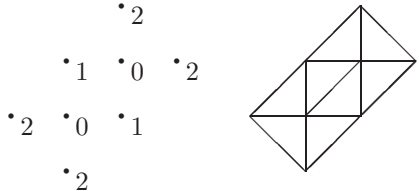
\begin{center}
\setlength{\unitlength}{3039sp}%
\begingroup\makeatletter\ifx\SetFigFont\undefined%
\gdef\SetFigFont#1#2#3#4#5{%
  \reset@font\fontsize{#1}{#2pt}%
  \fontfamily{#3}\fontseries{#4}\fontshape{#5}%
  \selectfont}%
\fi\endgroup%
\begin{picture}(3345,1548)(868,-1276)
{\thinlines
\put(1351,-1111){\circle*{50}}
}%
{\put(1801,-661){\circle*{50}}
}%
{\put(901,-661){\circle*{50}}
}%
{\put(1351,-211){\circle*{50}}
}%
{\put(1801,-211){\circle*{50}}
}%
{\put(2251,-211){\circle*{50}}
}%
{\put(1801,239){\circle*{50}}
}%
\put(2851,-661){\line( 1,-1){450}}
\put(3301,-1111){\line( 1, 1){900}}
\put(4201,-211){\line(-1, 1){450}}
\put(3751,239){\line(-1,-1){900}}
\put(3751,239){\line( 0,-1){900}}
\put(3751,-661){\line(-1, 0){900}}
\put(2851,-661){\line( 1, 1){900}}
\put(3301,-1111){\line( 0, 1){900}}
\put(3301,-211){\line( 1, 0){900}}
\put(4201,-211){\line(-1,-1){900}}
\put(3751,-211){\line(-1,-1){450}}

\put(1426,-1261){\makebox(0,0)[lb]{$2$}}
\put(976,-811){\makebox(0,0)[lb]{$2$}}
\put(1426,-811){\makebox(0,0)[lb]{$0$}}
\put(1876,-811){\makebox(0,0)[lb]{$1$}}
\put(1426,-361){\makebox(0,0)[lb]{$1$}}
\put(1876,-361){\makebox(0,0)[lb]{$0$}}
\put(2326,-361){\makebox(0,0)[lb]{$2$}}
\put(1876, 89){\makebox(0,0)[lb]{$2$}}
{\put(1351,-661){\circle*{50}}
}%
\end{picture}%
\end{center}
Then, $M$ is a genus $2$ curve with $6$ points removed, decomposed into $8$ pairs-of-pants. The Legendre dual yields the following decomposition:
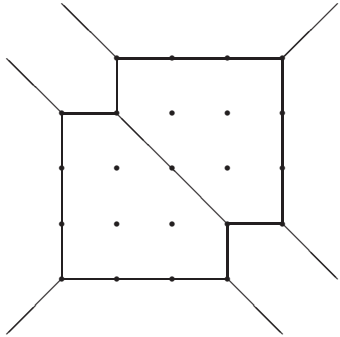
\begin{center}
\setlength{\unitlength}{3039sp}%
\begin{picture}(2724,2724)(-461,-2023)
\thinlines
{\put(  1,-211){\line(-1, 1){450}}
}%
{\put(1351,-1111){\circle*{50}}
}%
{\put(1801,-661){\circle*{50}}
}%
{\put(901,-661){\circle*{50}}
}%
{\put(1351,-211){\circle*{50}}
}%
{\put(1801,-211){\circle*{50}}
}%
{\put(1801,239){\circle*{50}}
}%
{\put(1351,239){\circle*{50}}
}%
{\put(901,239){\circle*{50}}
}%
{\put(451,239){\circle*{50}}
}%
{\put(451,-211){\circle*{50}}
}%
{\put(901,-211){\circle*{50}}
}%
{\put(  1,-211){\circle*{50}}
}%
{\put(451,-661){\circle*{50}}
}%
{\put(  1,-661){\circle*{50}}
}%
{\put(  1,-1111){\circle*{50}}
}%
{\put(451,-1111){\circle*{50}}
}%
{\put(901,-1111){\circle*{50}}
}%
{\put(1801,-1111){\circle*{50}}
}%
{\put(  1,-1561){\circle*{50}}
}%
{\put(451,-1561){\circle*{50}}
}%
{\put(901,-1561){\circle*{50}}
}%
{\put(1351,-1561){\circle*{50}}
}%
{\put(2251,689){\line(-1,-1){450}}
\put(1801,239){\line( 0,-1){1350}}
\put(1801,-1111){\line(-1, 0){450}}
\put(1351,-1111){\line( 0,-1){450}}
\put(1351,-1561){\line(-1, 0){1350}}
\put(  1,-1561){\line( 0, 1){1350}}
\put(  1,-211){\line( 1, 0){450}}
\put(451,-211){\line( 0, 1){450}}
\put(451,239){\line( 1, 0){1350}}
}%
{\put(1351,-1111){\line(-1, 1){900}}
}%
{\put(  1,-1561){\line(-1,-1){450}}
}%
{\put(1351,-1561){\line( 1,-1){450}}
}%
{\put(1801,-1111){\line( 1,-1){450}}
}%
{\put(451,239){\line(-1, 1){450}}
}%
{\put(1351,-661){\circle*{50}}
}%
\end{picture}%
\end{center}
Hence, the fibre $X = W^{-1}(0) \subset \XX$ has eight irreducible components, of which the two compact ones are isomorphic to $\C P^1 \times \C P^1$ blown up in one point. This can be considered as a degenerate version of the kind of pictures drawn in \cite{kapustin-katzarkov-orlov-yotov09}, and it also appears naturally when one applies a version of the SYZ construction to obtain a mirror for $(\C^*)^2 \times \C$ blown up at $M \times \{0\}$ \cite{abouzaid-auroux-katzarkov10}.
\end{example}

We now turn to categorical aspects. $M$ is an exact symplectic manifold with a canonical (up to homotopy) trivialization of its canonical bundle. Let $\F(M)$ be its Fukaya category. Here, we restrict the class of objects to compact, exact, graded, and {\it Spin} Lagrangian submanifolds, so that $\F(M)$ is a $\Z$-graded $A_\infty$-category defined over $\C$. On the other hand, there is a natural triangulated category $D^{\C^*}(W)$ associated to $W$, the category of equivariant Landau-Ginzburg branes \cite[Section 3]{orlov06}. By using dg quotients \`a la \cite{drinfeld02}, one can also define an underlying differential graded category $\Der^{\C^*}(W)$.

\begin{conjecture} \label{th:pp-1}
There is a cohomologically full and faithful embedding $\F(M) \hookrightarrow \Der^{\C^*}(W)$.
\end{conjecture}

The algebro-geometric object $(\XX,W)$ on the right hand side is glued together from pieces which follow an essentially unique standard model. In a suitable abstract sense, the same can be said of the associated category $\Der^{\C^*}(W)$. One could therefore imagine an approach to Conjecture \ref{th:pp-1} where the Fukaya category is viewed in the same way, and glued together from standard pieces corresponding to pair-of-pants in our decomposition. A first step towards constructing restriction functors which would govern the gluing process was taken in \cite{abouzaid-seidel07}. For a partial analysis of the Fukaya category of a single pair-of-pants, see \cite{seidel08b} (for $n = 1$) and \cite{sheridan10b} (for all $n$).

\begin{remark}
Let $\PP$ be the fan in $\R^{n+1}$ consisting of the cones $\PP_k = \bigcup_{t \in \R^+} \{t\} \times t P_k$. The associated toric variety is canonically isomorphic to $\XX$. Projection to $\R^+$ yields a toric map $\XX \rightarrow \mathbb{A}^1$, which is the same as our previously defined $W$. This description (for which see \cite[Section 3]{nishinou-siebert06} or \cite[Section 1.2]{gross-siebert08}) shows that $W: \XX \rightarrow \C$ depends only on the decomposition $P = \bigcup_k P_k$, and not on the function $\phi$ inducing that decomposition (the additional information provided by $\phi$ amounts to fixing an ample line bundle on $\XX$). Correspondingly, a construction of $M$ which uses only the $P_k$ is given in \cite[Corollary 5.1 and Remark 5.2]{mikhalkin04} (in this framework, $\phi$ determines an automorphism of $M$, namely the monodromy map obtained by moving the parameter $\epsilon$ around a small circle).
\end{remark}

\begin{remark}
One can change the $\C^*$-action on $\XX$ by composing it with any one-parameter subgroup of the natural torus which leaves $W$ invariant (which gives roughly a $\Z^n$ worth of choices), and such a modification would affect $\Der^{\C^*}(W)$. There are corresponding variations of $\F(M)$, obtained by changing the homotopy class of the trivialization of the canonical bundle.
\end{remark}

\begin{remark} \label{th:monodromy}
The embedding from Conjecture \ref{th:pp-1} can't be an equivalence, for two different reasons. First, $\Der^{\C^*}(W)$ is triangulated, and secondly (which is more important from a geometric perspective), it contains objects whose endomorphism spaces are infinite-dimensional. To improve the situation, one can enlarge $\F(M)$ to the so-called wrapped Fukaya category \cite{fukaya-seidel-smith07b, abouzaid-seidel07}, which includes some noncompact Lagrangian submanifolds, and then pass to its triangulated closure. It is a possible, if somewhat hazardous, notion that the outcome might be actually equivalent to $\Der^{\C^*}(W)$.
\end{remark}

\subsection*{Hypersurfaces in abelian varieties}
There is a variant of the previous picture in which one considers infinite periodic triangulations.
The objects encountered in this way are a little more interesting. In particular, on the symplectic geometry side, one gets closed (but nonsimply connected) manifolds.

Our starting point, in the viewpoint of \cite[Section 5]{mikhalkin-zharkov06}, is a polarized tropical abelian variety. Concretely, this will be given by specifying an integer lattice $\Gamma \subset \Z^n$ and a positive definite integral symmetric matrix $g$ of size $n$, such that $\gamma \cdot g\gamma \in 2\Z$ for $\gamma \in \Gamma$. In addition, suppose that we have a function $\phi_\Z: \Z^n \rightarrow \Z$ such that $\phi_\Z(v) - \half v \cdot gv$ is $\Gamma$-periodic; equivalently,
\begin{equation} \label{eq:quasip}
\phi_\Z(v + \gamma) = \phi_\Z(v) + v \cdot g\gamma + \half \gamma \cdot g\gamma
\end{equation}
for $v \in \Z^n$, $\gamma \in \Gamma$. Consider the entire function on $\C^n$ defined, for some $\epsilon>0$, by
\begin{equation} \label{eq:theta}
F(z) = \sum_{v \in \Z^n} \epsilon^{-2\pi i v \cdot z + 2\pi \phi_\Z(v)}.
\end{equation}
Because it's a Fourier series, this is $\log(\epsilon)^{-1}\Z^n$-periodic, and by \eqref{eq:quasip} it also satisfies
\begin{equation} \label{eq:f-auto}
F(z+ig\gamma) = F(z)\epsilon^{2\pi i z \cdot \gamma - \pi \gamma \cdot g\gamma}.
\end{equation}
Hence, its zero set descends to a hypersurface $M$ in the torus $T = \C^n/\log(\epsilon)^{-1}\Z^n \oplus ig(\Gamma)$. If we equip $M$ with the constant K{\"a}hler form
\begin{equation}
\omega = \log(\epsilon)\, d\mathrm{re}(z) \cdot g^{-1}\, d\mathrm{im}(z)
\end{equation}
then the homology class of $M$ is Poincar{\'e} dual to $[\omega]$ (since \eqref{eq:theta} is a linear combination of theta-functions for the factor of automorphy \eqref{eq:f-auto}, see for instance \cite[Chapter 2]{birkenhage-lange}). We will now impose further conditions. Namely, there should be a $\Gamma$-periodic tiling of $\R^n$ by integer simplices $P_k$ with volume $1/n!$, with the following property. Let $\phi: \R^n \rightarrow \R$ be the unique extension of $\phi_\Z$ to a function which is continuous, and affine on each $P_k$. We ask that $\phi$ must be convex, and not differentiable at any point lying on two different simplices. Then, $M$ will be smooth for sufficiently small values of $\epsilon$.

Take the Legendre transform of $\phi$, but where the standard scalar product in \eqref{eq:legendre} is replaced by $g$. The resulting function $\psi$ is convex, takes on integer values at integer points, and has the same periodicity property as $\phi$. We associate to it a graded ring $\hat\RR$, by modifying the previous construction for polytopes in the following way. Consider points $v \in \frac{1}{k}\Z^{n+1}$ such that $v_{n+1} \geq \psi(v_1,\dots,v_n)$, modulo the equivalence relation
\begin{equation} \label{eq:sim}
v \sim \big(v_1 + \gamma_1,\dots,v_n + \gamma_n,v_{n+1} + (v_1,\dots,v_n) \cdot g\gamma + \half \gamma \cdot g\gamma)
\end{equation}
for $\gamma \in \Gamma$. Each such point has a height, which is $v_{n+1} - \psi(v_1,\dots,v_n)$, and this is unchanged under \eqref{eq:sim}. When defining $\hat\RR^k$, we allow infinite formal linear combinations of the $v$, as long as their heights go to infinity. The result is that $\hat\RR^k$ is a finite rank free module over $\C\lb t \rb$, where $t$ acts as before by raising the last coordinate by $1/k$. Given two equivalence classes which give rise to generators of $\hat\RR$, we fix a point $v$ representing the first one, then allow all possible representatives $w$ of the second one under \eqref{eq:sim}: the (infinite) sum of the resulting terms \eqref{eq:product} defines the product of our generators. By taking an appropriate version of $\Proj$, we can associate to $\hat\RR$ a formal scheme $\hat\XX$, which comes with a proper flat morphism $\hat{W}: \hat\XX \rightarrow \hat{\mathbb{A}}^1$ to the formal disc. This construction has a long history in algebraic geometry, most of which is frankly outside the competence of this author (the starting point is \cite{mumford72}; for an expository account of recent developments, see \cite{brion05}). Because of our assumptions, $\hat\XX$ itself is smooth, and its zero fibre has toric components, which have normal crossings with each other.

\begin{example} \label{th:node}
Take $\Gamma = \Z$, $g = 2$, and $\phi_\Z(v) = v^2$. After substituting $\epsilon = \exp(\tau/2i)$ and $z = 2i \zeta/\tau$, one sees that $F$ is precisely the Jacobi theta function for $\C/\Z \oplus \tau\Z$, hence $M$ will consist of a single point. The Legendre transform yields $\psi = \phi$. The quotient $R = \hat\RR/t\hat\RR$ has generators $a \in R^1$, $b \in R^2$ and $c \in R^3$ corresponding to points $(0,0)$, $(\frac{1}{2},\frac{1}{2})$ and $(\frac{1}{3},\frac{1}{3})$, respectively, and the relation
\begin{equation}
\textstyle abc = b^3 + c^2 = (\frac{1}{2},\frac{1}{2}) + (\frac{1}{3},\frac{1}{3}) \in R^6.
\end{equation}
Therefore, the zero-fibre $W^{-1}(0) = \Proj(R)$ of our degeneration is a rational curve with a node (lying at $b = c = 0$ in the coordinates introduced above).
\end{example}

\begin{example}
Take $\Gamma = 2\Z \times 2\Z \subset \Z^2$, $g = 4 \cdot \mathrm{Id}$, and
\begin{equation}
\phi_\Z(v) = \begin{cases}
2\|v\|^2 & \text{if $(v_1,v_2) = (0,0)$ mod $2$}, \\
2\|v\|^2 - 2 & \text{if $(v_1,v_2) = (0,1)$ or $(1,0)$ mod $2$}, \\
2\|v\|^2 - 3 & \text{if $(v_1,v_2) = (1,1)$ mod $2$.}
\end{cases}
\end{equation}
This decomposes $\R^2/\Gamma$ into eight triangles
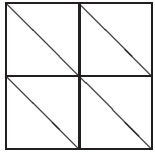
\begin{center}
\setlength{\unitlength}{4000sp}%
\begin{picture}(924,924)(1189,-523)
\thinlines
{\put(1201,-511){\framebox(900,900){}}
}%
{\put(1651,389){\line( 0,-1){900}}
}%
{\put(1201,-61){\line( 1, 0){900}}
}%
{\put(1201,-61){\line( 1,-1){450}}
}%
{\put(1201,389){\line( 1,-1){900}}
}%
{\put(1651,389){\line( 1,-1){450}}
}%
\end{picture}%
\end{center}
In this case, $M$ is a genus five curve, correspondingly decomposed into eight pairs-of-pants. The central fibre of the degeneration associated to the Legendre dual $\psi$ has four irreducible components, each of which is a copy of $\C P^2$ blown up at three points.
\end{example}

By construction, all our hypersurfaces $M$ have ample canonical bundle. More precisely, that bundle carries a connection whose curvature is proportional to $\omega$. This allows one to define the Fukaya category $\F(M)$ as a $\Z/2$-graded category over $\C$, by appropriately restricting the class of Lagrangian submanifolds under consideration (compare \cite{seidel08b}; for the more familiar monotone case which has the opposite sign of $c_1(M)$, see \cite{oh93,wehrheim-woodward06}). On the other side, the main difference with respect to the previous situation is that there is no longer a $\C^*$-action on $\hat\XX$ covering that on the base. Hence, we necessarily have to work with the ordinary differential graded category of Landau-Ginzburg branes $\Der(\hat{W})$, which is $\Z/2$-graded.

\begin{conjecture} \label{th:pp-2}
There is a cohomologically full and faithful embedding $\F(M) \hookrightarrow \Der^\pi(\hat{W})$, where the superscript $\pi$ stands for split-closure.
\end{conjecture}

\begin{remark}
The split-closure (also called idempotent or Karoubi completion) is necessary even in simple situations like Example \ref{th:node}. In that case, $\Der^\pi(\hat{W})$ is generated by a single object, the structure sheaf at the singular point of the zero-fibre. The endomorphisms of that object form a two-dimensional Clifford algebra. Hence, in the split-closure, the object splits into two summands which are isomorphic up to a shift. Each summand corresponds to the point object in the Fukaya category $\F(M)$ (this problem would not have been present in the parallel case of Conjecture \ref{th:pp-1}). More generally, passing to the idempotent completion is also a great computational simplification, since we know split-generators of $\F(M)$ in many cases \cite{seidel04}.
\end{remark}

\begin{remark}
If we set $q = \log(\epsilon)$, then the torus $T \iso \C^n/\Z^n \oplus iq g(\Gamma)$ can be defined for complex values of $q$ (in the left half-plane), and is invariant under $q \mapsto q+i$. The same holds for the hypersurface $M$, which therefore comes with a monodromy map (compare Remark \ref{th:monodromy}).
The same conclusion can be reached under slightly weaker assumptions than the ones used before: instead of asking for $g$ to be integral, it is sufficient that $g(\Gamma) \subset \Z^n$. However, the Legendre transformed function $\psi$ then no longer necessarily takes integer values on integers, so the construction of the mirror would have to be modified. 
\end{remark}

\section{The analytic topology\label{sec:analytic}}

\subsection*{The local picture for torus fibrations}
We start by reviewing some elementary material about the symplectic geometry of Lagrangian torus fibrations, see for instance \cite{duistermaat80}. Consider $M = T^n \times \R^n$ with its standard symplectic form. This is a Lagrangian torus fibration over $B = \R^n$, whose symplectic geometry is tied closely to the $\Z$-affine geometry of the base. For instance, take a function $H$ on $B$, pull back it to $M$, and consider the associated Hamiltonian system. This system is one-periodic if and only if $H$ is $\Z$-affine, which means of the form $H(b) = t + v \cdot b$ for some $t \in \R$ and $v \in \Z^n$. Similarly, consider $\Z$-affine automorphisms of $B$, which are maps of the form $\phi(b) = u + Ab$, where $u \in \R^n$ and $A \in GL_n(\Z)$. Denote the group of such automorphisms by $G$. They then naturally lift to symplectic automorphisms of $M$. The last-mentioned observation has a more interesting converse:

\begin{theorem}[Benci \protect{[unpublished]} for the convex case, Sikorav \cite{sikorav91} in general] Let $U,V \subset \R^n$ be connected open subsets with trivial first Betti number. If $T^n \times U$ is symplectically isomorphic to $T^n \times V$, there is a $\phi \in G$ such that $\phi(U) = V$.
\end{theorem}

Let $\Lambda_q$ be the Novikov field in one variable $q$, and with coefficients in $\C$. This is the field of formal series $\sum_r a_r q^r$, where $r \in \R$ and $a_r \in \C$, such that for each $C \in \R$, there are only finitely many $r \leq C$ with $a_r \neq 0$. It comes with an obvious real-valued valuation $\val_q: \Lambda_q \rightarrow \R \cup \{+\infty\}$. Kontsevich and Soibelman \cite{kontsevich-soibelman00} introduced the following sheaf $\Olo_B$ of commutative $\Lambda_q$-algebras on $B$.  For any connected open $U$, $\Olo_B(U)$ consists of series $f = \sum_{v \in \Z^n} c_v z_1^{v_1} \cdots z_n^{v_n}$, $c_v \in \Lambda_q$, such that the following property holds for each $b \in U$:
\begin{equation} \label{eq:valuation}
\text{\it given $C \in \R$, there are only finitely many $v \in \Z^n$ such that
$\val_q(c_v) - b \cdot v \leq C$.}
\end{equation}
More canonically, the exponents $v$ naturally belong to the lattice of integer cotangent vectors on $B$. This shows that the pairing in \eqref{eq:valuation} is well-defined without reference to any basis.

\begin{example}
Take $n = 1$ and $U = (-\infty,0)$. For a function $f = \sum_v c_v z^v$ to lie in $\Olo_B(U)$, the sequence $\val_q(c_v)$ must satisfy $\lim_{v \rightarrow -\infty} \val_q(c_v)/v \rightarrow -\infty$ (super-linear decay) and $\liminf_{v \rightarrow +\infty} \val_q(c_v)/v \geq 0$ (sub-linear growth). This is the nonarchimedean analogue of the notion of holomorphic functions on the punctured open unit disc in $\C$.
\end{example}

Algebras of this kind are common in nonarchimedean geometry, but they can also be related to  symplectic homology (references for symplectic homology are \cite{viterbo97a,cieliebak-floer-hofer95}, or for instance the more recent \cite{cieliebak-frauenfelder-oancea09}). Namely, let $H$ be a function on $B$, and $X_H$ the Hamiltonian vector field of its pullback to $M$, which is parallel to the torus fibres. Take the zero-section $L = \{0\} \times \R^n \subset M$. We are interested in length one chords, which are trajectories $x: [0,1] \rightarrow M$ of $X_H$ with both endpoints on $L$. There is precisely one such chord $x_b$ in the fibre over every point $b$ where $dH$ is an integral cotangent vector. The action functional for such chords is given by
\begin{equation} \label{eq:action}
A_H(x_b) = H(b) - dH_b \cdot b.
\end{equation}
Assume that all $x_b$ are nondegenerate, which means that the Hessians $D^2 H_b$ are nonsingular.
Then they have well-defined Maslov indices, which are just the Morse indices of $H_b$. One defines the Floer cochain space $\CF^*(L,L;H)$ as follows. Its elements are formal series $\sum_{b,r} a_{b,r} q^r x_b$, where $a_{b,r} \in \C$, the $x_b$ are chords of a given Maslov index, $r$ runs over the real numbers, and the following condition holds: for each $C \in \R$, there are only finitely many nonzero $a_{b,r}$ such that $A_H(x_b) + r \leq C$. If there are only finitely many $x_b$ overall, then $\CF^*(L,L;H)$ is just the $\Lambda_q$-vector space generated by them. However, in general it is somewhat bigger, being complete for the $q$-adic topology.

Now take an open subset $U \subset B$ which is bounded and strongly convex. The last-mentioned condition means that it is convex, with smooth boundary, and that the Gauss map $N: \partial U \rightarrow S^{n-1}$ is a diffeomorphism. One can then find a continuous function $\bar{H}: \bar{U} \rightarrow \R$ which is constant on the boundary, restricts to a smooth and strictly convex function $H = \bar{H}|U$ on the interior, and such that the derivatives grow to infinity in normal direction near the boundary. More precisely, we want
\begin{equation} \label{eq:radial-growth}
\begin{aligned}
& \lim_{b \rightarrow c} \|dH_{b}\| = \infty, \\
& \lim_{b \rightarrow c} dH_{b}/\|dH_{b}\| = N(c)
\end{aligned}
\end{equation}
uniformly for $c \in \partial U$ (a general reference for the construction of convex functions is \cite{blocki97}, but our specific situation is actually much easier than the one considered in that paper). This in particular implies that $dH: U \rightarrow \R^n$ is a diffeomorphism. Temporarily ignoring the fact that the function is defined only partially, we can then consider the Floer cochain space $\CF^*(L,L;H)$, which in this case is concentrated in degree zero. The chords $x_b$ are naturally parametrized by $v = dH_b \in \Z^n$, so that we can write elements of $\CF^0(L,L;H)$ as series $\sum_{v,r} a_{v,r} q^r z^v$. Due to \eqref{eq:action} and the fact that $H$ is bounded, the growth condition says that for each $C$, there are only finitely many nonzero $a_{v,r}$ such that $r - v \cdot b \leq C$. This is just the specialization of \eqref{eq:valuation} to points $b$ where the chords appear. However, because of \eqref{eq:radial-growth} such points cluster everywhere near the boundary, which implies that
\begin{equation} \label{eq:o-f}
\CF^0(L,L;H) \iso \Olo_B(U).
\end{equation}
Hamiltonian functions such as $H$ are tricky from the point of view of Floer cohomology, so one may have trouble using this approach directly (which would be necessary if one wanted to refine \eqref{eq:o-f} to an isomorphism of rings rather than vector spaces, for instance). However, there are known workarounds where one first considers a certain class of functions which are constant outside $U$, and then passes to the limit, see again \cite{cieliebak-floer-hofer95, cieliebak-frauenfelder-oancea09}.

\subsection*{Global torus fibrations}

Again following \cite{kontsevich-soibelman00}, we now globalize the picture above. Namely, let $B$ be a $\Z$-affine manifold. This structure can either be described through its sheaf $\Aff_B$ of $\Z$-affine functions (seen as a subsheaf of the continuous function sheaf $\Cont_B$), or more concretely by specifying an atlas with transition functions in $G$. There is a canonical symplectic manifold $M$ associated to $B$, defined as the quotient of $T^*B$ by its fibrewise integer lattice $T^*_\Z B$, and a natural sheaf of $\Lambda_q$-algebras $\Olo_B$ in the base, locally isomorphic to the one introduced above for $B = \R^n$.

From now on, we assume that $B$ is closed. Let $\Der(\Olo_B)$ be the derived dg category of $\Olo_B$-module sheaves. This is a dg category over $\Lambda_q$, whose underlying cohomological category is the standard unbounded derived category, and is unique up to quasi-equivalence. On the symplectic side, consider the Fukaya category $\F(M)$, which is an $A_\infty$-category over $\Lambda_q$. To precisely specify this category, we have to make at least two additional remarks. First of all, for any compatible almost complex structure on $M$, the canonical bundle $\Kan_M = \Lambda_\C^{\mathrm{top}} TM$ can be identified with the complexified pullback of $\Lambda_\R^{\mathrm{top}} TB$. In particular, $\Kan_M^2$ has a canonical trivialization up to homotopy, which we can use to define the gradings in the Fukaya category. Secondly, we can pull back the Stiefel-Whitney class $w_2(B)$ to $M$, and use Lagrangian submanifolds which are {\it Pin} relative to this class \cite{fooo} to fix the signs. With this in mind, we have:

\begin{theorem}[Kontsevich-Soibelman \cite{kontsevich-soibelman00}]
The full subcategory of $\F(M)$ consisting of Lagrangian sections of our torus fibration embeds in a cohomologically full and faithful way into $\Der(\Olo_B)$.
\end{theorem}

\begin{conjecture}[Kontsevich-Soibelman \cite{kontsevich-soibelman00}] \label{th:ks-conj}
All of $\F(M)$ embeds into $\Der(\Olo_B)$ in the same sense.
\end{conjecture}

The conjecture follows from the theorem above if one can prove that the Lagrangian sections are split-generators for the Fukaya category. This is known when $B$ is a circle \cite{polishchuk-zaslow98} or a two-dimensional square torus \cite{abouzaid-smith09}.

We also need a slight generalization. A {\em signed $\Z$-affine manifold} is a $\Z$-affine manifold $B$ together with a torsor $L^\sigma$ over the local system $T^*_{(1/2) \Z}B / T^*_{\Z} B \iso T^*_{\Z}B \otimes \Z/2$ (equivalently, one can express this additional datum as a lift of the transition maps from $G$ to a finite extension $G^\sigma = (\Z/2)^n \rtimes G$). This naturally gives rise to a Lagrangian torus fibration $M \rightarrow B$ together with an antisymplectic fibrewise involution, whose fixed point set can be identified with $L^\sigma$. On the other hand, we can associate to our signed $\Z$-affine structure a sheaf of matrix algebras $\Olo_B^\sigma$ of size $2^n$ over $\Olo_B$. For $B = \R^n$ with the trivial torsor, rows and columns of the matrices are labeled by elements of $(\Z/2)^n$, and the $(i,j)$ entries are elements of $z_1^{(i_1-j_1)/2}\cdots z_n^{(i_n-j_n)/2} \Olo_B$. In general, rows and columns are labeled by points in the fibre of $L^\sigma$. The analogue of Conjecture \ref{th:ks-conj} says that $\F(M)$ should embed into the derived category of $\Olo_B^\sigma$-module sheaves. For the trivial torsor, this reduces to the original statement by way of Morita equivalence.

In the remainder of the paper, we will consider only the case when $B$ is one-dimensional. In that case, $T^*_{(1/2) \Z}B / T^*_{\Z} B$ is canonically isomorphic to the constant sheaf $\Z/2$, so $L^\sigma$ is just a double cover of $B$. We want to tweak the definition of $\Olo_B^\sigma$ slightly, by turning it into a sheaf of $\Z/2$-graded algebras, where the off-diagonal pieces $z^{1/2} \Olo_B$ are odd, and the others even. The covering transformation of $L^\sigma$ acts by permuting the rows as well as the columns of $\Olo_B^\sigma$.

\subsection*{The pair-of-pants}
Our next step is to adapt the mirror of the pair-of-pants from Section \ref{sec:zariski} to the nonarchimedean framework. Consider tropical affine three-space, which is $B = (\R^+)^3$ with the sheaf $\Olo_B$ of functions $f = \sum_{v \in \Z^3} c_v z^v$ satisfying the same condition as in \eqref{eq:valuation}, and with the following additional holomorphicity requirement: if an open connected subset $U \subset B$ contains a point $b$ with $b_k = 0$ for some $k \in \{1,2,3\}$, then all $f \in \Olo_B(U)$ have $c_v = 0$ for all $v$ such that $v_k < 0$. In particular, global sections consist of formal power series with rapidly decreasing $\Lambda_q$-coefficients. Take the superpotential $W = -z_1 z_2 z_3 \in \Olo_B(B)$, and its three obvious matrix factorizations, which are
\begin{equation}
\xymatrix{
E_1 = \{
\cdots \Olo_B \ar[rr]^-{z_1} && \Olo_B \ar[rr]^-{-z_2 z_3} && \Olo_B \cdots\}
}
\end{equation}
and its cyclic permutations $E_2$, $E_3$. The endomorphisms of $E = E_1 \oplus E_2 \oplus E_3$ form a sheaf of $\Z/2$-graded differential algebras $\A_B = \End(E)$. If we forget the differential and the grading, then this is the matrix algebra over $\Olo_B$ of rank $6$.

Let $U_3 = (\R^+)^2 \times (0,\infty) \subset B$ be the complement of the third coordinate plane, $I_3 = (0,\infty)$ the open half-line, and $\iota_3: I_3 \rightarrow U_3$ the embedding given by $\iota_3(b) = (0,0,b)$. Consider the trivial double cover of $I_3$ and the associated sheaf of $\Z/2$-graded algebras $\Olo^\sigma_{I_3}$.

\begin{lemma} \label{th:sheaf-quasi}
There is a quasi-isomorphism of sheaves of differential algebras,
\begin{equation} \label{eq:phi-push}
\phi_3: (\iota_3)_*\Olo^\sigma_{I_3} \longrightarrow \A_B|U_3.
\end{equation}
\end{lemma}

In informal language, the idea behind this statement is that on the region where $z_3 \neq 0$, the superpotential becomes a nondegenerate quadratic function of $z_1,z_2$. Matrix factorizations of such functions form a semisimple category, see for instance \cite{kapustin-li03}, which is described by the left hand side of \eqref{eq:phi-push}. More concretely, one first uses $z_3^{-1}$ to build a contraction of $E_3|U_3$, which shows that $\A_B|U_3$ is quasi-isomorphic to its subalgebra of endomorphisms of $E_1 \oplus E_2$. One then writes down an explicit quasi-isomorphism
\begin{equation}
(\iota_3)_*\Olo^\sigma_{I_3} \longrightarrow \End(E_1 \oplus E_2)|U_3,
\end{equation}
as follows. Elements on the left hand side are two-by-two matrices $(f_{ij}(z_3))$, where the diagonal entries are series with $z_3^n$ exponents, while the off-diagonal entries have $z_3^{n+1/2}$ exponents. We take a function $f_{kk}(z_3)$ to the corresponding diagonal degree zero endomorphism of $E_k$, respectively. The off-diagonal entries go to the odd homomorphisms between $E_1$ and $E_2$, given respectively by
\begin{equation} \label{eq:1-to-2}
\xymatrix{
E_1|U = \{\cdots \Olo_B \ar[drr]^{\;\;\;\;z_3^{-1/2} f_{12}(z_3)} \ar[rr]^-{z_1} && \ar[drr]^{\quad z_3^{1/2} f_{12}(z_3)} \Olo_B \ar[rr]^-{-z_2 z_3} && \Olo_B \cdots\} \\
E_2|U = \{\cdots \Olo_B \ar[rr]_-{z_2} && \Olo_B \ar[rr]_-{-z_3 z_1} && \Olo_B  \cdots\}
}
\end{equation}
and the same formula in the other direction. Consider the transposition $\tau_{12}$ which exchanges the first two coordinates. Then we have a commutative diagram
\begin{equation} \label{eq:swap-ends}
\xymatrix{
\Olo^\sigma_{I_3} \ar[d] \ar[rr]^-{\phi_3} && \A_B|U_3 \ar[d]^{\tau_{12}^*} \\
\Olo^\sigma_{I_3} \ar[rr]^-{\tau_{12}^*(\phi_3)} && \tau_{12}^*\A_B|U_3,
}
\end{equation}
where the left hand $\downarrow$ is the involution of $\Olo^\sigma_{I_3}$ associated to the nontrivial covering transformation.

The {\em tropical pair-of-pants} is
\begin{equation} \label{eq:trivalent}
T = (\R^+ \times \{0\}) \cup (\{0\} \times \R^+) \cup \{(-b,-b) \;:\; b \geq 0\} \subset \R^2.
\end{equation}
Take the embedding $\iota: T \rightarrow B$ given by $(b,0) \mapsto (b,0,0)$, $(0,b) \mapsto (0,b,0)$, $(-b,-b) \mapsto (0,0,b)$, and set
\begin{equation}
\Olo_T^\sigma = \iota^*\A_B.
\end{equation}
In general, as a pullback sheaf this would be defined as a direct limit, but Lemma \ref{th:sheaf-quasi} shows that in this case all the maps in the direct system are quasi-isomorphisms. Consider the boundary of a thickening of $T \subset \R^2$, together with its projection to $T$, which is a double cover except for the preimage of the vertex:
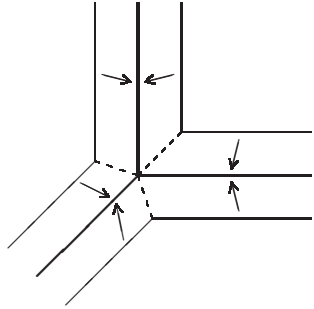
\begin{equation} \label{eq:retract}
\setlength{\unitlength}{0.00050000in}
\begingroup\makeatletter\ifx\SetFigFont\undefined%
\gdef\SetFigFont#1#2#3#4#5{%
  \reset@font\fontsize{#1}{#2pt}%
  \fontfamily{#3}\fontseries{#4}\fontshape{#5}%
  \selectfont}%
\fi\endgroup%
{\renewcommand{\dashlinestretch}{30}
\begin{picture}(3184,3199)(0,-10)
\drawline(1212,687)(1137,1062)
\drawline(1198.777,982.573)(1137.000,1062.000)(1110.524,964.923)
\thicklines
\drawline(1362,1362)(312,312)
\thinlines
\drawline(3162,1812)(1812,1812)(1812,3162)
\drawline(912,3162)(912,1512)(12,612)
\drawline(612,12)(1512,912)(3162,912)
\dashline{60.000}(1812,1812)(1362,1362)
\dashline{60.000}(1362,1362)(1512,912)
\dashline{60.000}(1362,1362)(912,1512)
\drawline(987,2412)(1287,2337)
\drawline(1188.773,2315.172)(1287.000,2337.000)(1210.601,2402.485)
\drawline(1737,2412)(1437,2337)
\drawline(1513.399,2402.485)(1437.000,2337.000)(1535.227,2315.172)
\drawline(2412,1737)(2337,1437)
\drawline(2315.172,1535.227)(2337.000,1437.000)(2402.485,1513.399)
\drawline(2412,987)(2337,1287)
\drawline(2402.485,1210.601)(2337.000,1287.000)(2315.172,1188.773)
\drawline(762,1287)(1062,1137)
\drawline(961.377,1137.000)(1062.000,1137.000)(1001.626,1217.498)
\thicklines
\drawline(3162,1362)(1362,1362)(1362,3162)
\end{picture}
}
\end{equation}
A cyclic ordering of the edges of $T$, or equivalently an orientation of the plane into which it is embedded, determines a trivialization of the double cover over $I_3$, and simultaneously a quasi-isomorphism $\Olo_{I_3}^\sigma \htp \Olo_T^\sigma|I_3$.

\subsection*{Metrized graphs}
Let $C$ be a graph with trivalent vertices, but possibly with some semi-infinite edges. A {\em signed tropical structure} on $C$ is given by a $\Z$-affine structure on each edge, together with a map $L^\sigma \rightarrow C$ from a topological one-manifold $L^\sigma$ which is a double cover away from the vertices, and has the same local structure near those vertices as in \eqref{eq:retract}. This restricts to give a signed $\Z$-affine structures on the edges, and it also allows us to glue together the resulting cylinders to form a symplectic surface $M$, in an essentially unique way (which means up to symplectic automorphisms which are Hamiltonian isotopic to the identity). The Fukaya category $\F(M)$ is a $\Z/2$-graded $A_\infty$-category over $\Lambda_q$ (if $M$ is not a torus, there are smaller subcategories defined over $\C$, versions of which we used previously in Conjectures \ref{th:pp-1} and \ref{th:pp-2}; however, those do not fit in naturally with the present argument).

On the other hand, one can associate to $C$ a sheaf of differential $\Z/2$-graded algebras $\Olo_C^\sigma$, by using the standard matrix algebra sheaf over the edges, the model $\Olo_T^\sigma$ near each vertex, and the quasi-isomorphisms from Lemma \ref{th:sheaf-quasi} to glue the two together. The general theory of sheaves of dga's can be quite sophisticated, because one wants to allow homotopies between the gluing maps (compare for instance \cite{toen-vezzosi04}), but the present case is much more straightforward, because of the particularly simple topology of the base. There is an associated dg category $\Der(\Olo_C^\sigma)$ of sheaves of dg modules (appropriately derived), and our main idea is expressed by the following:

\begin{conjecture} \label{th:c}
$\F(M)$ embeds cohomologically fully and faithfully into $\Der(\Olo_C^\sigma)$.
\end{conjecture}

\begin{remark}
It is helpful to compare this with previously proposed descriptions of the Fukaya category of closed higher genus surfaces. Several such constructions have appeared in the literature, starting with the pioneering \cite{kapustin-katzarkov-orlov-yotov09} (and indeed, we've given another one in Section \ref{sec:zariski}, for surfaces appearing as theta-divisors). In all such cases, the mirror is a Landau-Ginzburg model whose singular set is an algebraic curve consisting of rational components whose intersection models the pair-of-pants decomposition of the surface. However, the local algebro-geometric structure near that set varies, and even though one expects all of these descriptions to be equivalent (after idempotent completion), it is by no means obvious that this is true. One possible advantage of Conjecture \ref{th:c} is that the model is constructed in a more hands-on way, with less extraneous information.
\end{remark}


\end{document}